\documentclass[12pt]{article}
\usepackage{amsmath,graphicx,amssymb,amscd}
\usepackage[latin1]{inputenc}
\begin{document}

\title{\textbf{Pairs of subsets of spheres and Cartesian products thereof with the same distribution of distance}}
\author{
Ricardo Garc\'\i a-Pelayo\thanks{%
E-mail: r.garcia-pelayo@upm.es} \\
\\
ETS de Ingenier\'ia Aeron\'autica \\
Plaza del Cardenal Cisneros, 3 \\
Universidad Polit\'ecnica de Madrid \\
Madrid 28040, Spain\\}
\date{}
\maketitle

\begin{abstract}
We prove the following three statements: 1) Let $(A, \overline A)$ be a partition of the spherical surface $S^n$ into two measurable sets. Let $st_A$ and $st_{\overline A}$ be their measure density functions of distance. Then $|st_A - st_{\overline A}|$ depends only on the difference of their $n$-areas. 2) If the spherical surface $S^n$ is divided in two measurable subsets $A$ and $\overline A$ of equal $n$-surface, then these two subsets have the same distribution of distance. 3) Let there be a pair $(S, S')$ of subsets of a sphere $S^{n}$ such that $st_S = st_{S'}$. Then their complementary subsets satisfy $st_{\overline S} = st_{\overline S'}$ and $st_{S, \overline S} = st_{S', \overline S'}$, where $st_{A, B}$ is the measure density function of distance between a point in $A$ and a point in $B$. Furthermore, it is shown that the statements remain true when $S^n$ is substituted by the Cartesian product $S^{n_1} \times ... \times S^{n_r}$ endowed with the metric which is naturally inherited from its factors.

\end{abstract}

\bigskip

Some pairs of non-congruent polygons are known to share some geometric properties (see e. g. Figures 1.4 and 3.8 of \cite{Gardner1995}). In particular, non-congruent polygons \cite{Mallows1970, RGPN2016a} and polyhedra \cite{RGPN2016a} exist such that their distributions of distance or of length of chords are the same. It was recently conjectured by this author (conjectures 6.1, 6.2 and 6.3 of \cite{RGPN2016a}) that the extension of some results on non congruent subsets of regular polyhedra with the same distribution of distance to non congruent subsets of spheres of arbitrary dimension were true. We show in this article that these conjectures are right except that the probability density functions of distance mentioned in the conjectures are actually measure density functions of distance.

In the propositions and proofs to follow the distance used is the Euclidean distance of the embedding space, $R^{n+1} \supset S^n$. This we do for technical and graphical convenience. All the propositions of this article are propositions about measure density functions of distance being equal. All these propositions to follow remain true if the Euclidean distance is substituted by the angular distance, that is, the angle $\alpha \in [0, \pi]$ subtended by the shortest arc of maximum circle which joins two given points. The reason for that is that $\alpha$ and the Euclidean distance $d_E$ determine each other through the relation $2 \sin (\alpha/2) = d_E$ ($S^n$ has unit radius), and two measure density functions of the Euclidean distance are equal iff the two measure density functions of the angular distance are equal.

In \cite{RGPN2016a} the collection of distances was defined. While this definition makes sense for distances between finite sets, it cannot be extended straightforwardly to infinite sets. For any kind of set, the general idea is that each distance $\ell$ appears a number of times proportional to the number of couples of points which are a distance $\ell$ apart. In particular, the measure of any Lebesgue measurable set $A$ induces a Stieltjes measure (\cite{Kolmogorov1960}, Vol. 2, p. 13) on the interval $[0, {\rm{diam}}(A)]$. This cumulative Stieltjes measure at $\ell$ is $St(\ell) = \int_A d^n r_1\ \int_A d^n r_2\ H(\ell - d_E(\vec r_1, \vec r_2))$, where $d_E$ is the Euclidean distance and $H$ is the Heaviside function, which is 1 for positive arguments, and 0 otherwise. It reaches $\mu(A)^2$ when $\ell = {\rm{diam}}(A)$. Its derivative is the measure density function $st(\ell) = \int_A d^n r_1\ \int_A d^n r_2\ \delta(\ell - d_E(\vec r_1, \vec r_2))$, where $\delta$ is Dirac's delta function.

{\bf{Notation.}} $St_A$ is the cumulative measure function of distance between points of a set $A$. $st_A$ is the measure density function of distance between points of a set $A$. $St_{A,B}$ is the cumulative measure function of distance between a point in $A$ and a point in $B$. $st_{A,B}$ is the measure density function of distance between a point in $A$ and a point in $B$. A ball $BB(\vec r, r) \subset S^n$ of radius $r$ and center $\vec r \in S^n$ is the intersection of a ball in $B(\vec r, r) \subset R^{n+1}$ and $S^n$, that is, $BB(\vec r, r) \equiv B(\vec r, r) \cap S^n$.

{\bf{Lemma.}} The swapping of two balls between a measurable subset of $S^n$ and its complementary does not change the difference of their measure density functions of distance.

\begin{figure}[ht]
  \begin{center}
  \includegraphics[width=0.6 \textwidth]{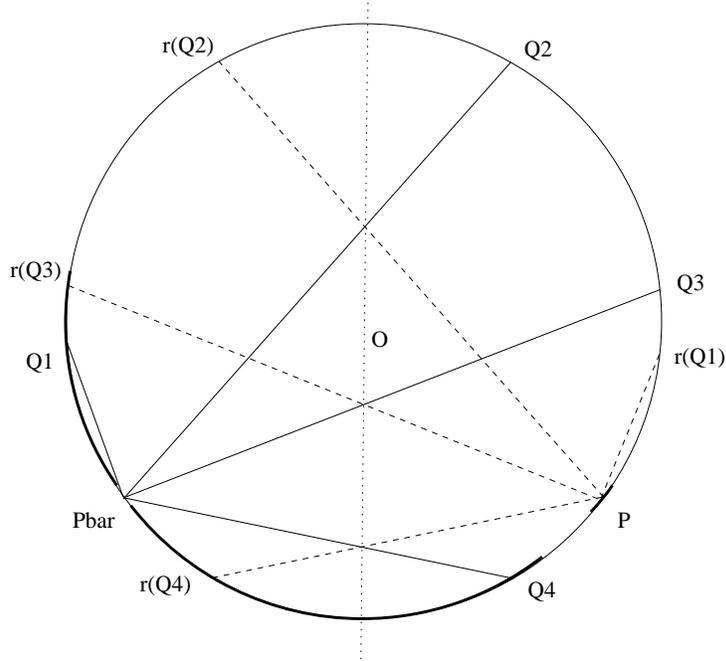}\caption{The swapping of two balls between $A$ and $\bar A$ does not change the difference $\rho_A - \rho_{\overline A}$. The $n$-hyperplane which bisects the angle $(P, O, \overline P)$ is represented by a dotted line. The 4 possible cases are depicted.}
  \end{center}
\end{figure}

{\bf{Proof.}} Take balls $BB({\overline P}, r) \subset {\overline A}$ and  $BB(P, r) \subset A$ and swap them. Consider the angle $(P, O, \overline P)$, where $O$ is the centre of the sphere. There is a unique $n$-hyperplane which bisects the angle $(P, O, \overline P)$. Any point $Q$ on $S^n$ has an image $r(Q)$ when reflected with respect to this $n$-hyperplane. For any point $\overline P' \in BB({\overline P}, r)$, its reflection $P' \in BB(P, r)$ is such that the segments $\overline P' Q$ and $P' r(Q)$ have the same length.

In the swapping of the balls $BB({\overline P}, r)$ and  $BB(P, r)$ the distances of the lengths of the form $\overline P' Q$ are either added to the collection of distances of $A$ ($\overline P' Q_2$ and $\overline P' Q_3$ in Fig. 1), or are subtracted from the collection of distances of $\overline A$ ($\overline P' Q_1$ and $\overline P' Q_4$ in Fig. 1). In either case they contribute positively to $st_A - st_{\overline A}$. Likewise, the  distances of the lengths of the form $P' r(Q)$ contribute negatively to $st_A - st_{\overline A}$. Since the segments of the forms  $\overline P' Q$ and  $P' r(Q)$ are each other's reflection,  $st_A - st_{\overline A}$ remains unchanged. This may be expressed in mathematical notation as follows:

$$
\big( St_A - St_{\overline A} \big) (\ell) = \int_{BB(\overline P, r)} d^n \overline P'\ \int_{S^n} d^n Q\
\Big( H \big[ \ell - \overline P' Q \big] - H \big[ \ell - P' r(Q) \big] \Big) = 0.
$$

{\bf{QED}}

{\bf{Main Lemma.}} Let $(A, \overline A)$ be a partition of the spherical surface $S^n$ into two measurable sets.
Let $st_A$ and $st_{\overline A}$ be their measure density functions of distance. Then $|st_A - st_{\overline A}|$ depends only on the difference in their $n$-areas.

{\bf{Proof.}} We have to prove that for any two partitions of $S^n$ into two measurable subsets, say $(A, \overline A)$ and $(B, \overline B)$, if $\mu(B) = \mu(A)$, then $St_A - St_{\overline A} = St_B - St_{\overline B}$.

Pick a pair of balls $BB(\vec r_{A1}, r_1)$ and $BB(\vec r_{B1}, r_1)$, of equal radius $r_1$, such that $BB(\vec r_{A1}, r_1)\\ \subset A-B$ and $BB(\vec r_{B1}, r_1) \subset B-A$ and such that $r_1$ is as large as possible. Swap the balls. Then choose a pair of balls $BB(\vec r_{A2}, r_2)$ and $BB(\vec r_{B2}, r_2)$ such that $BB(\vec r_{A2}, r_2) \subset A-B-BB(\vec r_{A1}, r_1)$ and $BB(\vec r_{B2}, r_2) \subset B-A-BB(\vec r_{B1}, r_1)$ and $r_2$ is as large as possible. Swap the balls and repeat again.

If this procedure reaches an end, the lemma is proven. If not, the measure \\ $\mu(\sqcup_i^\infty BB(\vec r_{Ai}, r_i)) = \mu(\sqcup_j^\infty BB(\vec r_{Bj}, r_j))$ covered by the pairs must reach a limit $\mu'$. If $\mu' < \mu(A-B) = \mu(B-A)$, then the sets $A - B - \sqcup_i^\infty BB(\vec r_{Ai}, r_i)$ and $B - A - \sqcup_j^\infty BB(\vec r_{Bj}, r_j)$ are measurable sets of non zero measure. Then each of them contains a ball of positive radius, and thus a further couple of balls of non zero measure could be added to $\{ BB(\vec r_{Ai}, r_i) \}_i^\infty$ and $\{ BB(\vec r_{Bj}, r_j) \}_j^\infty$, thus it must be $\mu' = \mu(A-B) = \mu(B-A)$. But if $\mu' = \mu(A-B) = \mu(B-A)$, the lemma is proven. {\bf{QED}}

{\bf{Theorem 1.}} If the spherical surface $S^n$ is divided in two measurable subsets $A$ and $\overline A$ of equal $n$-surface, then these two subsets have the same distribution of distance.

{\bf{Proof.}} We know that there are complementary subsets of equal measure which are also congruent, for example the two hemispheres $H$ and $\overline H$, and therefore have the same distribution of distance. Since $\mu(A) - \mu(\overline A) = \mu(H) - \mu(\overline H)$, the theorem follows from the main lemma. {\bf{QED}}

{\bf{Theorem 2.}} Let there be a pair $(S, S')$ of subsets of a sphere $S^{n}$ such that $st_S = st_{S'}$. Then their complementary subsets satisfy $st_{\overline S} = st_{\overline S'}$ and $st_{S, \overline S} = st_{S', \overline S'}$.

{\bf{Proof.}} Since the measure density functions of distance $st_S$ and $st_{S'}$ are the same, it must be $\mu(S) = \mu(S')$. Then
$\mu(S) + \mu(\overline S) = \mu(S') + \mu(\overline S') \Rightarrow \mu(\overline S) = \mu(\overline S')$ and the first statement follows from the lemma.

As for the second statement, $st_{S^n} = st_{S} + st_{\overline S} + st_{S, \overline S} =
st_{S'} + st_{\overline S'} + st_{S', \overline S'}$, and the second statement follows. {\bf{QED}}

{\bf{Theorem 3.}} The Cartesian product $S^{n_1} \times ... \times S^{n_r}$ can be endowed with a metric $d$ as follows

$$
d\big( (\vec r_1, ..., \vec r_r), (\vec r'_1, ..., \vec r'_r) \big) = f\big( d_E(\vec r_1, \vec r'_1), ..., d_E(\vec r_r, \vec r'_r) \big),
$$
where $f$ is a function such that the axioms of metric are satisfied. If $f$ is measurable, then the main lemma and theorems 1 and 2 are true when $S^n$ is substituted by $S^{n_1} \times ... \times S^{n_r}$.

{\bf{Proof.}} If $f$ is measurable, then the measure density function of distance for any measurable subset of $S^{n_1} \times ... \times S^{n_r}$ is defined. What makes the main lemma work is that for any pair of points in $S^n$ there is a reflection of $S^n$ which maps the points of the pair onto each other. This is also true for any pair of points in $S^{n_1} \times ... \times S^{n_r}$, because there is such a reflection for the coordinates of each factor of the Cartesian product. {\bf{QED}}

Examples of Theorem 3 are

$$
d\big( (\vec r_1, ..., \vec r_r), (\vec r'_1, ..., \vec r'_r) \big) = \sqrt{d_E(\vec r_1, \vec r'_1)^2 + ... + d_E(\vec r_r, \vec r'_r)^2}
$$
or, for $S \times S$, the flat Clifford torus \cite{Kholodenko2013} (p. 118) endowed with the distance

$$
d\big( (\theta_1, \theta_2), (\theta'_1, \theta'_2) \big) = \sqrt{ \alpha(\theta_1, \theta'_1)^2 + \alpha(\theta_2, \theta'_2)^2},
$$
where $\alpha$ is the angular distance mentioned above:

$$
\alpha(\theta, \theta') = \min \big( |\theta - \theta'|, 2 \pi - |\theta - \theta'| \big).
$$

The methods used in section 5 of \cite{RGPN2016a} to construct pairs of non congruent convex bodies with the same distribution of distance cannot be applied, at least straightforwardly, to any pair of measurable subsets of $S^{n_1} \times ... \times S^{n_r}$.


\begin{thebibliography}{1}

\bibitem{Gardner1995}
Richard~J. Gardner.
\newblock {\em {Geometric Tomography}}.
\newblock Cambridge University Press, 1995.

\bibitem{Mallows1970}
C.~L. Mallows and J.~M.~C. Clark.
\newblock {Linear-Intercept Distributions Do Not Characterize Plane Sets}.
\newblock {\em Journal of Applied Probability}, 7(1):240--244, 1970.

\bibitem{RGPN2016a}
\rm{Ricardo Garc\'ia-Pelayo}.
\newblock {Pairs of Subsets of Regular Polyhedra with the Same Distribution of
  Distance}.
\newblock {\em Applied Mathematical Sciences}, 10(26):1285--1297, 2016. \\
  http://dx.doi.org/10.12988/ams.2016.6396.

\bibitem{Kolmogorov1960}
A.~N. Kolmogorov and S.~V. Fomin.
\newblock {\em Elements of the Theory of Functions and Functional Analysis}.
\newblock Dover, 1999.

\bibitem{Kholodenko2013}
Arkady~L. Kholodenko.
\newblock {\em {Applications of Contact Geometry and Topology in Physics}}.
\newblock World Scientific, 2013. ISBN:978-981-4412-08-7.

\end{thebibliography}

\end{document}